\documentclass[12pt]{article}
\usepackage{amsmath}
\usepackage{amscd}
\usepackage{graphicx}
\usepackage{amssymb,amsfonts,amsthm,amscd,mathrsfs}
\usepackage{verbatim}
\usepackage[T2A]{fontenc}
\usepackage[utf8]{inputenc}
\usepackage[russian]{babel}
\usepackage[dvipsnames]{xcolor} 
\usepackage{bm}
\usepackage{authblk}

\usepackage{amsmath,amsthm,amssymb}
\usepackage{amssymb,amsmath,amsthm}
\newtheorem{theorem}{Теорема}
\numberwithin{theorem}{section}

\newtheorem{remark}{Замечание}

\numberwithin{lemma}{section}
\numberwithin{proposition}{section}
\numberwithin{corollary}{section}
\numberwithin{definition}{section}
\numberwithin{example}{section}
\numberwithin{exercise}{section}
\numberwithin{problem}{section}
\numberwithin{remark}{section}

\begin{document}
\title{
Testing of symmetry of innovations in autoregression 
}

\author{M.~V.~Boldin, A.~R.~Shabakaeva    
\footnote{Moscow State Lomonosov Univ., Dept. of Mech. and Math., Moscow, Russia\\
e-mail: boldin$_{-}$m@hotmail.com, ~ 
shabakaevalmira@gmail.com}}

\date{ }
\maketitle

\textbf{Abstract}

We consider a stationary linear $AR(p)$ model with zero mean. The  autoregression parameters as well as the distribution function (d.f.) $G(x)$ of innovations are unknown. We consider two situations.

In the first situation the observations are a sample from a stationary solution of $AR(p)$. Interesting  and essential problem is to test symmetry of $G(x)$ with respect to zero. If hypothesis of symmetry is valid then it is possible to construct nonparametric estimators of $AR(p)$ parameters, for example, GM-estimators, minimum distance estimators and others. 

First of all  we estimate unknown parameters of autoregression and find residuals. Based on them we construct a kind of empirical d.f., which is a counterpart of empirical d.f of the unobservable innovations. Our test statistic is the functional of  omega-square type from this residual empirical d.f. Its asymptotic d.f. under the hypothesis and the local alternatives are found.

In the second situation the observations subject to gross errors (outliers). The distribution of outliers is unknown, their intensity is $O(n^{-1/2})$, $n$ is the sample size. We test the symmetry of innovations again but by constructing the Pearson's type statistic. Its asymptotic d.f. under the hypothesis and the local alternatives are found. We establish the asymptotic robustness of this test as well.



\section{Введение и постановка задачи}\leavevmode
Рассмотрим линейную AR(p) модель 
$$
u_t = \beta_1 u_{t-1} + \dots + \beta_p u_{t-p} + \varepsilon_t, ~ ~ t \in \mathbb{Z}. \eqno(1.1)
$$

В (1.1) инновации $\{ \varepsilon_t \}$ - независимые одинаково распределенные случайные величины (н.о.р.сл.в.) с неизвестной функцией распределения (ф.р.) $G(x)$; $E \varepsilon_1 = 0$, $0 < E \varepsilon_1^2 < \infty$; $\bm\beta = (\beta_1, \dots, \beta_p)^T$ -- вектор неизвестных параметров, для которых соответствующее характеристическое уравнение имеет корни, лежащие в единичном круге;
размерность модели $p$ предполагается известной. \\
Эти требования считаются выполненными всегда и далее особо не оговариваются.

Мы будем рассматривать две ситуации. В первой наблюдения $u_{1-p}, u_{2-p}, \dots, u_n$ суть выборка из стационарного решения (1.1). По этим наблюдениям хотим проверить гипотезу
$$
\mathbb{S} ~ : ~ \varepsilon_1 \overset{d}{=} -\varepsilon_1,  \,\,\text{т.е.}\,\, G(x) \,\, \text{симметрична относительно нуля. } \eqno(1.2)
$$
Для проверки $\mathbb{S}$ мы построим тест типа омега-квадрат со статистикой $\hat{\omega}_n^2,$ найдем асимптотическое при $n\to\infty$ распределение $\hat{\omega}_n^2$ при гипотезе и локальных альтернативах. Чтобы описать эти альтернативы,  будем предполагать, что ф.р. $G(x)$  зависит от числа наблюдений $n$ и представляется смесью симметричной относительно нуля ф.р. $P(x)$ и несимметричной ф.р. $Q(x)$ :
$$
G(x) = A_n(x) := (1 - \rho_n) P(x) + \rho_n Q(x). \eqno(1.3)
$$
В (1.3) ф.р. $P(x)$ и $Q(x)$ неизвестны,  $\rho_n = \min{\{1, \rho \, n^{-1/2}\}}, $ $\rho \geq 0,$ $\rho$ неизвестно.\\
Если выполнено (1.3), будем говорить, что верна гипотеза $\mathbb{A}_n(\rho)$. При $\rho > 0$ $\mathbb{A}_n(\rho)$ будем понимать как локальную альтернативу к $\mathbb{S}$, ф.р. $A_n(x)$ при $\rho>0$ несимметрична. При $\rho = 0,$
разумеется, $\mathbb{A}_n(0)$ и $\mathbb{S}$ совпадают. В в этом случае (т.е. при гтпотезе $\mathbb{S}$) будем писать $P(x)$ вместо $G(x)$. Асимптотические распределения $\hat{\omega}_n^2$ будут найдены при $\mathbb{A}_n(\rho)$ сразу для всех $\rho \geq 0,$ что позволяет изучить тест одновременно и при гипотезе, и при локальных альтернативах. 

Отметим, что симметрия инноваций -- одно из основных предположений, позволяющих строить непараметрические ценки для $\bm\beta$. Например, $GM-$оценки, оценки минимального расстояния, знаковые и т.д., см. \cite{Koul}. Кроме того, для симметричной $G(x)$ можно строить симметризованные оценки $G(x)$ и основывать на них симметризованные тесты для гипотез о виде $G(x),$ см. \cite{Boldin2020}. Строить такие тесты, как и тесты для $\mathbb{S},$ нетривиальная и содержательная задача, поскольку инновации $\{\varepsilon_t\}$ ненаблюдаемы.\\
Соответствующие первой ситуации результаты представлены в Разделе 2. 

Во второй ситуации наблюдения за авторегрессией содержат грубые ошибки. А именно, наблюдаются 
$$y_t = u_t + z_t^{\gamma_n} \xi_t, ~ t = 1-p, \dots, n. \eqno(1.4)$$
В (1.4) $\{z_t^{\gamma_n}\}$ -- н.о.р. сл.в., распределенные по закону Бернулли, $z_1^{\gamma_n} \sim Br(\gamma_n)$, $\gamma_n = \min\{1, \gamma n^{-1/2}\}, ~ \gamma \geq 0, \gamma$ неизвестно;
$\{\xi_t\}$ -- н.о.р. сл.в. с неизвестным и произвольным распределением $\Pi$; $\{u_t\}$ -- выборка из стационарного решения (1.1); последовательности $\{u_t\}$, $\{z_t^{\gamma_n}\}$, $\{\xi_n\}$ независимы.
Последовательность $\{\xi_t\}$ интерпретируется как последовательность грубых ошибок (засорений). Схема $(1.4)$ -- локальный вариант известной схемы засорений для временных рядов, см. \cite{Martin}.

Задачей ставится опять проверка гипотезы $\mathbb{S}$ из (1.2) по наблюдениям $\{y_t\}$. Мы построим тест типа хи-квадрат для $\mathbb{S}$, изучим асимптотическое поведение тестовой статистики при $\mathbb{A}_n(\rho)$ из $(1.3)$, установим асимптотическую качественную робастность теста. Эти результаты приведем в Разделе 3.  

\section{Проверка симметриии в схеме без засорений}
Пусть наблюдения $u_{1-p}, \dots, u_n$ -- выборка из стационарного решения уравнения $(1.1)$. Построим по этим наблюдениям тест для проверки гипотезы $\mathbb{S}$ из $(1.2)$.

Пусть $\hat{\bm\beta}_n = (\hat{\beta}_{1 n}, \dots, \hat{\beta}_{p n})^T$ -- оценка вектора $\beta$, требования к ней  уточним далее. Оценками ненаблюдаемых инноваций $\varepsilon_1, \dots, \varepsilon_n$ возьмем остатки 
$$
\hat{\varepsilon}_t := u_t - \hat{\beta}_{1 n} u_{t-1} - \dots - \hat{\beta}_{p n} u_{t-p}, ~ ~ t = 1, \dots, n;
$$
оценкой ф.р. $G(x)$ -- остаточную эмпирическую ф.р. 
$$
\hat G_n(x) = n^{-1} \sum\limits_{t = 1}^n I(\hat{\varepsilon}_t \leq x), ~ ~ x \in \mathbb{R}^1,\,\,I(\cdot) -- \text{индикатор события.}
$$
Тестоввой статистикой для $\mathbb{S}$ возьмем
$$\hat{\omega}_n^2 = n \int_{-\infty}^{\infty} [\hat{G}_n(x) + \hat{G}_n(-x) - 1]^2 d \hat{G}_n(x). $$
Вычислять $\hat{\omega}_n^2$ удобно по формуле 
$$\hat{\omega}_n^2 = \sum\limits_{t=1}^n \left[\hat{G}_n(-\hat{\varepsilon}_{(t)}) - \dfrac{n-t+1}{n} \right]^2, $$
где $\hat{\varepsilon}_{(1)} \leq \dots \leq \hat{\varepsilon}_{(n)}$ -- упорядоченные остатки.

Статистика $\hat{\omega}_n^2$ -- аналог статистики $\omega_n^2$ для проверки симметрии, см. \cite{Orlov1972}, которую можно было бы построить по самим $\varepsilon_1, \dots, \varepsilon_n$. А именно,
$$
\omega_n^2 = n \int_{-\infty}^{\infty} [G_n(x) + G_n(-x) - 1]^2 d G_n(x),
$$
где $G_n(x) = n^{-1} \sum\limits_{t=1}^n I(\varepsilon_t \leq x)$ -- эмпирическая ф.р. $\varepsilon_1, \dots, \varepsilon_n$.

Пусть $v(t), ~ t \in [0, 1],$ -- броуновский мост, т.е. гауссовский процесс с нулевым средним и ковариацией $\min{\{t, s\}} - t s$. Для $\mathbb{A}_n(\rho)$ из (1.3) положим 
$$\delta(t): = \rho \, [ Q\,(P^{-1}\,(t)\,) - t ],$$
где $P^{-1}(\cdot)$ -- обратная к $P(\cdot)$ функция. Из результатов \cite{Orlov1972} следует:\\
\hangindent = 0.5 cm
    \textit{при гипотезе $\mathbb{S}$ для непрерывной $G(x)$ 
    $$\omega_n^2 \xrightarrow[]{~d~} \int_0^1 [v(t) + v(1-t)]^2 dt, ~ ~ n \to \infty,$$
    а при альтернативе $\mathbb{A}_n(\rho)$ с $\rho > 0$ для непрерывных $A_n(x)$ и $\delta(t)$
    $$\omega_n^2 \xrightarrow[]{~d~} \int_0^1 [v(t) + v(1-t) + \delta(t) + \delta(1-t)]^2 dt, ~ ~ n \to \infty.$$ }

Нам понадобятся следующие условия.\\
\hangindent = 0.5 cm
    $\textbf{Условие } \mathbf{(i)}. ~~ $ Ф.р. $P(x)$ и $Q(x)$ имеют нулевые средние и конечные дисперсии.\\    
    $\textbf{Условие } \mathbf{(ii)}. ~$ Ф.р. $P(x)$ и $Q(x)$ дважды дифференцируемы с ограниченными вторыми производными.\\    
    $\textbf{Условие } \mathbf{(iii)}.$ $\hat{\bm\beta}_n$ -- такая оценка $\bm\beta$, что для фиксированного 
     $\rho \geq 0$ при $\mathbb{A}_n(\rho)$ 
    $$
    n^{1/2} (\hat{\bm\beta}_n -\bm \beta) = O_P(1), ~ ~ n \to \infty. \eqno(2.1)
    $$
Если выполнено Условие (i), то оценка наименьших квадратов $\hat{\bm\beta}_{n, LS}$ , построенная по $u_{1-p}, \dots, u_n,$ удовлетворяет $(2.1),$ см.\cite{Boldin2019}. Более того, 
$$
n^{\frac{1}{2}} (\hat{\bm\beta}_{n, LS} - \bm\beta) \xrightarrow[]{~d~} \mathcal{N}(0, E(\varepsilon_1^0)^2 \, \mathcal{K}^{-1}), ~ ~ n \to \infty. \eqno(2.2)
$$
В $(2.2)$ $E(\varepsilon_1^0)^2$ -- дисперсия ф.р. $P(x),$ $\mathcal{K} = (k_{i j}),$ $k_{i j} = E (u_0^0 \, u_{i - j}^0), ~ i, j = 1, \dots, p,$\\
$\{u_t^0\}$ -- стационарное решение $(1.1)$ с инновациями $\{\varepsilon_t^0\},$ имеющими ф.р. $P(x)$.

Далее нам понадобится результат из \cite{Boldin2019}:

\hangindent=0.5cm
\textit{если выполнены Условия (i) - (iii), то при $\mathbb{A}_n(\rho)$ с любым $\rho \geq 0$}
$$
\sup_{x \in \mathbb{R}^1} |n^{1/2} [\hat{G}_n(x) - G_n(x)]| = o_{P}(1), ~ n \to \infty. \eqno(2.3)
$$

В силу $(2.3)$ асимптотическое распределение $\hat{\omega}_n^2$ при $\mathbb{A}_n(\rho),\,\rho \geq 0$, совпадает с асимптотическим распределением $\omega_n^2$, т.е. верна 
\begin{theorem}
$ 1^0$. Пусть верна гипотеза    $\mathbb{S}$. Пусть для ф.р. $P(x)$ выполнены Условия (i) -- (ii) и Условие (iii) c $\rho=0$ для $\hat{\bm\beta}_{n}$. Тогда    
    $$
    \hat{\omega}_n^2 \xrightarrow[]{~d~} \int_0^1 [v(t) + v(1-t)]^2 d t, ~ n \to \infty.
    $$

$2^0$. Пусть при некотором $\rho>0$ верна альтернатива $\mathbb{A}_n(\rho)$. Пусть выполнены Условия (i) - (iii), и фукция 
 $\delta(t), ~ t \in [0, 1]$, непрерывна. Тогда 
   $$
    \hat{\omega}_n^2 \xrightarrow[]{~d~} \int_0^1 [v(t) + v(1-t) + \delta(t) + \delta(1-t)]^2 d t, ~ n \to \infty.
    $$ 
    
\end{theorem}

Для $0 < \alpha < 1$ критическое множество возьмем в виде 
$$\hat{\omega}_n^2 > c_{1 - \alpha},  \eqno(2.4)$$
где $c_{1 - \alpha}$ -- $(1 - \alpha)$-квантиль предельной ф.р. $\omega_n^2$.
Эта предельная ф.р. табулирована в \cite{Mart1978}. Асимптотический уровень теста $(2.4)$ равен $\alpha$.

\begin{remark}
    При фиксированной альтернативе к гипотезе $\mathbb{S}$  ф.р. инноваций несимметрична относительно нуля. Пусть  она, кроме того, имеет нулевое среднее, конечную дисперсию и ограниченную вторую производную. Пусть $n^{1/2} (\hat{\bm\beta}_n -\bm \beta) = O_P(1), ~ ~ n \to \infty$.- Тогда статистика $\hat{\omega}_n^2 $ расходится по вероятности к бесконечности. Это означает, что тест $(2.4)$ состоятелен против такой фиксированной альтернативы.
\end{remark}

\begin{remark}
    Помимо $\hat{\omega}_n^2$ можно строить аналоги и других статистик, упомянутых в \cite{Orlov1972}. Например,
    $$
    \hat{D}_n = \sup_{x \in \mathbb{R}^1} |n^{1/2} [\hat{G}_n(x) - \hat{G}_n(-x) - 1]|.
    $$
    Асимптотическое распределение $\hat{D}_n$ при $\mathbb{A}_n(\rho),\,\rho \geq 0$, такое же, как у статистики 
    $$
    D_n = \sup_{x \in \mathbb{R}^1} |n^{\frac{1}{2}} [G_n(x) + G_n(-x) -1]|.
    $$
\end{remark}

\section{Проверка симметрии в схеме с засорениями}

Пусть наблюдаются $y_{1-p}, \dots, y_n$, определенные в $(1.4)$.
Построим по этим наблюдениям тест для проверки гипотезы $\mathbb{S}$ из $(1.2)$.

Прежде всего построим по $\{y_t\}$ оценку $\hat{\bm\beta}_n^y = (\hat{\beta}_{1 n}^y, \dots, \hat{\beta}_{p n}^y)^T$ вектора $\bm\beta$. Далее $R \geq 0, ~ \Gamma \geq 0$ -- произвольные конечные числа.\\
\hangindent = 0.5 cm
\textbf{Условие } $\mathbf{(iv)}$. $\hat{\bm\beta}_n^y$ -- такая оценка $\bm\beta$, для которой при $\mathbb{A}_n(\rho)$ равномерно по\\ $\rho \leq R, ~ \gamma \leq \Gamma$ выполнено 
$$
n^{1/2} (\hat{\bm\beta}_n^y - \bm\beta) = O_P(1), ~ n \to \infty. \eqno(3.1)
$$
Если $E \varepsilon_1^4 < \infty$, то соотношение $(3.1)$ выполнено, например, для о.н.к. $\hat{\bm\beta}_{n, LS}^y,$ поскольку
$$
n^{1/2} (\hat{\bm\beta}_{n, LS}^y - \bm\beta) \xrightarrow[]{~d~} \mathcal{N}(-\gamma \, \mathcal{K}^{-1} \, \beta, E(\varepsilon_1^0)^2 \mathcal{K}^{-1}), ~ n \to \infty, \eqno(3.2)
$$
и сходимость по распределению равномерная по $\rho \leq R, ~ \gamma \leq \Gamma.$ 

Найдем остатки 
$\hat{\varepsilon}_t^y = y_t - \hat{\beta}_{1 n}^y \, y_{t-1} - \dots - \hat{\beta}_{p n}^y \, y_{t-p}, ~ ~ t = 1, \dots, n, $

и построим по ним остаточную эмпирическую ф.р.
$\hat{G}_n^y(x) = n^{-1}\sum\limits_{t = 1}^n I(\hat{\varepsilon}_t^y \leq x).$

Справедливо следующее утверждение ( \cite{Bold}, Теорема 2.1, Следствие 2.1):\\
\hangindent = 0.7 cm
{\it Если выполнены Условия (i), (ii), (iv), то при $\mathbb{A}_n(\rho),\,\rho\geq 0$, для фиксированного $x \in \mathbb{R}^1$
$$
n^{1/2} [\, \hat{G}_n^y(x) - G_n(x) \, ] - \gamma \Delta(x, \Pi) = o_P(1), ~ n \to \infty, \eqno(3.3)
$$
равномерно по $\rho \leq R, ~ \gamma \leq \Gamma.$ В $(3.3)$ 
$$
\Delta(x, \Pi) = \sum\limits_{j=0}^p [E P(x + \beta_j \xi_1) - P(x)], ~ ~ \beta_0 = -1.
$$
}

Разложение $(3.3)$ для $\hat{G}_n^y$ справеливо только при каждом фиксированном $x$, потому его недостаточно для исследования основанных на $\hat{G}_n^y(x)$ статистик типа $\hat{\omega}_n^2$ и $\hat{D}_n$ из Раздела 2. Однако достаточно для построения и исследования статистики типа хи-квадрат для проверки $\mathbb{S}$. Вот как она строится.

Введем полуинтервалы 
$$
B_j^+ = (x_{j-1}, x_j], ~ ~ j = 1, \dots, m, ~ m \geq 1, ~~~ 0 = x_0 < x_1 < \dots < x_m = \infty.
$$
Пусть $\{x_j\}$ таковы, что $p_j^+ := P(x_j) - P(x_{j-1}) > 0.$
Если ввести симметричные полуинтервалы $B_j^- = (-x_{j}, -x_{j-1}]$, то при гипотезе $\mathbb{S}$ для непрерывной $P(x)$
$p_j^- := P(-x_{j-1}) - P(-x_j) = p_j^+, ~ ~ j = 1, \dots, m.$

Пусть $\hat{\nu}_j^+ $ обозначает число остатков среди $\hat{\varepsilon}_1^y, \dots, \hat{\varepsilon}_n^y$, попавших в $B_j^+$, а $\hat{\nu}_j^-$ -- число остатков, попавших в $B_j^-.$
Интересующая нас тестовая статистика для $\mathbb{S}$ имеет вид
$$\hat{\chi}^2_n = \sum\limits_{j=1}^m \dfrac{(\hat{\nu}_j^+ - \hat{\nu}_j^-)^2}{2 \hat{\nu}_j^+}$$
Очевидно, 
$$
\dfrac{\hat{\nu}_j^+}{n} = \hat{G}_n^y(x_j) - \hat{G}_n^y(x_{j-1}), ~~ \dfrac{\hat{\nu}_j^-}{n} = \hat{G}_n^y(-x_{j-1}) - \hat{G}_n(-x_j).  \eqno(3.4)
$$

Пусть $\nu_j^{\pm}$ обозначает число инноваций $\varepsilon_1, \dots, \varepsilon_n$, попавших в $B_j^{\pm}.$
Поскольку 
$$
\hat{\chi}_n^2 = \sum\limits_{j=1}^m \{n^{1/2} \, (\frac{\hat{\nu}_j^+}{n} - \frac{\hat{\nu}^-}{n})\}^2/(2 \hat{p}_j^+), ~ ~ \hat{p}_j^+ = \frac{\hat{\nu}_j^+}{n},
$$
получаем в силу $(3.3)$ и $(3.4)$:
$$
\hat{\chi}_n^2 = \sum\limits_{j=1}^m \{n^{1/2} \, (\frac{{\nu}_j^+}{n} - \frac{{\nu}^-}{n}) + \delta_j(\Pi)\}^2/(2 \hat{p}_j^+)+o_P(1), ~ ~ n \to \infty. \eqno(3.5)
$$
$$\delta_j(\Pi): = \Delta(x_j, \Pi) - \Delta(x_{j-1}, \Pi) - [\Delta(-x_{j-1}, \Pi) - \Delta(-x_j, \Pi)].$$
В (3.5) $o_P(1)$ обозначает величину, стремящуюся по вероятности к нулю при $n \to \infty$ равномерно по $\gamma \leq \Gamma, ~ \rho \leq R.$ 

Соотношение $(3.5)$ сводит асимптотическое исследование статистики $\hat{\chi}_n^2$  к анализу главного члена в правой части (3.5).
Результат дается Теоремой 3.1. В ней 
$$
q_j^+ = Q(x_j) - Q(x_{j-1}), ~ ~ q_j^- = Q(-x_{j-1}) - Q(-x_j), ~ ~\bm  q^{\pm} = (q_1^{\pm}, \dots, q_m^{\pm})^T,
$$
$$
\delta(\bm\Pi) = (\delta_1(\Pi), \dots, \delta_m(\Pi))^T,
$$
$$\mathcal{P} \text{-- диагональная матрица,} ~~ \mathcal{P} = diag\{2p_1^+, \dots, 2p_m^+\}.$$
Через $F_k(x, \lambda^2)$ мы обозначаем ф.р. нецентрального хи-квадрат распределения с $k$ степенями свободы и параметром нецентральности $\lambda^2,$ а через $|\cdot|$ -- евклидову норму вектора.

\begin{theorem}
$ 1^0$. Пусть верна гипотеза    $\mathbb{S}$. Пусть для ф.р. $P(x)$ выполнены Условия (i) -- (ii) и Условие (iv) c $\rho=0$ для $\hat{\bm\beta}_n^y$. Тогда    
   $$
   \sup_{x \in \mathbb{R}^1, ~ \gamma \leq \Gamma} |P(\hat{\chi}_n^2 \leq x) - F_m(x, \lambda^2)| \to 0, ~ n \to \infty, \eqno(3.6)
   $$
    где $\lambda^2 = \gamma^2 \, |\, \mathcal{P}^{-\frac{1}{2}} \, \bm\delta(\Pi) \, |^2;$\\
    
$2^0$.Пусть верна альтернатива $\mathbb{A}_n(\rho),\,\,\rho>0$. Пусть выполнены Условия (i) - (ii), (iv). Тогда 
$$
\sup_{x \in \mathbb{R}^1, ~\gamma \leq \Gamma, ~\rho \leq R} |P(\hat{\chi}_n^2 \leq x) - F_m(x, \lambda^2)| \to 0, ~ n \to \infty, \eqno(3.7)
$$
    где $$
   \lambda^2 = |\mathcal{P}^{-1/2} [\rho(\bm q^+ - \bm q^-) + \gamma \bm\delta(\Pi)]|^2. \eqno(3.8)
    $$
\end{theorem}

Критическое множество для $\mathbb{S}$ возьмем в виде 
$$\hat{\chi}_n^2 > \chi_{1-\alpha}(m), \eqno(3.9)$$
где для $0 < \alpha < 1 ~~ \chi_{1 - \alpha}(m)$ -- $(1 - \alpha)$-квантиль ф.р. хи-квадрат с $m$ степенями свободы. 
В силу $(3.9)$  и Теоремы 3.1 асимптотическая мощность нашего теста есть
$$W(\rho, \gamma, \Pi, P, Q) = 1 - F_m(\chi_{1-\alpha}(m), \lambda^2),$$
где $\lambda^2 = \lambda^2(\rho, \gamma, \Pi, P, Q)$ задается формулой $(3.8).$ Известно, что
$$|F_k(x, \lambda_1^2) - F_k(x, \lambda_2^2)| \leq \sqrt{\dfrac{2}{\pi}}|\lambda_1 - \lambda_2|, \eqno(3.10),$$
так что из $(3.8)$, $(3.10)$ и определения вектора $\bm\delta(\Pi)$ получаем:
$$
\sup_{\Pi, \rho, Q} |W(\rho, \gamma, \Pi, P, Q) - W(\rho, 0, \Pi, P, Q)| \leq \sqrt{\dfrac{2}{\pi}} \, \gamma \, \sup_{\Pi} |\mathcal{P}^{-1/2}\bm\delta(\Pi)| \to 0, ~ \gamma \to 0. \eqno(3.11)
$$
Соотношение $(3.11)$ качественно означает, что при малых $\gamma$ равномерно по $\Pi$ (а также по $\rho$ и $Q$) асимптотические мощности в схемах с засорениями и без засорений близки. Это свойство означает асимптотическую качественную робастность теста.

Сделаем несколько замечаний.
\begin{remark}
    Если $\gamma = \rho = 0,$ то асимптотическое распределение для $\hat{\chi}_n^2$ есть обычное (центральное) распределение хи-квадрат с $m$ степенями свободы, $W(0, 0, \Pi, P, Q) = \alpha.$
    
    Если же распределение $\Pi$ симметрично, т.е. $\xi_1 \overset{d}{=} -\xi_1,$ то $\bm\delta(\Pi) = 0$, и асимптотическая мощность $W(\rho, \gamma, \Pi, P, Q)$ вовсе от $\Pi$ не зависит и совпадает с асимптотической мощностью теста в случае схемы без засорений: для $\rho \geq 0$ при симметричном распределении $\Pi$ и любых $\gamma\geq 0,\rho\geq 0$
    $$W(\rho, \gamma, \Pi, P, Q) = W(\rho, 0, \Pi, P, Q).$$
\end{remark}

\begin{remark}
    Рассмотрим фиксированную альтернативу для $\mathbb{S}$ c дополнительными предположениями как в Замечании 2.1. Если при каком-нибудь $j$ 
    $P(\varepsilon_1 \in B_j^+) \neq P(\varepsilon_1 \in B_j^-),$
    то $\hat{\chi}_n^2 \xrightarrow[]{P} \infty, ~ n \to \infty.$
    Т.е. тест (3.9)  состоятелен против такой фиксированной альтернативы.
\end{remark}

\begin{remark}
    Если вместо Условия $(iv)$ потребовать лишь, чтобы выполнялось $(3.1)$ при фиксированных $\rho$ и $\gamma$, то утверждение типа Теоремы 3.1 останется верным. Надо только в $(3.6)$ и $(3.7)$ заменить супремум по $x \in \mathbb{R}^1, ~\rho \leq R, ~\gamma \leq \Gamma$ на супремум по $x \in \mathbb{R}^1.$ 
    Например, о.н.к. $\hat{\beta}_{n, LS}^y$ удовлетворяет $(3.2)$ (а значит, и $(3.1)$) при фиксированных $\rho$ и $\gamma$ при минимальном условии $E \xi_1^2 < \infty.$ 
\end{remark}

\end{document}